\documentclass[12pt,reqno]{amsart}
\usepackage{amsmath, amsfonts, amssymb, amsthm, hyperref}
\usepackage{bm}
\allowdisplaybreaks[4]
\def\l{\left}
\def\r{\right}
\def\bg{\bigg}
\def\({\bg(}
\def\){\bg)}
\def\t{\text}
\def\f{\frac}
\def\i{\mathrm i}
\def\mo{{\rm{mod}\ }}

\def\sign{{\rm sign}}

\def\ls{\leqslant}

\def\bi{\binom}

\def\ve{\varepsilon}

\def\eq{\equiv}

\def\da{\delta}

\def\Proof{\noindent{\it Proof}}

\def\Z{\mathbb Z}

\def\Q{\mathbb Q}
\def\R{\mathbb R}

\def\p{\mathfrak p}
\def\1{{\bf 1}}

\def\pmod #1{\ ({\rm{mod}}\ #1)}

\def\<{\langle}
\def\>{\rangle}
\theoremstyle{plain}
\newtheorem{theorem}{Theorem}[section]
\newtheorem{lemma}{Lemma}[section]
\newtheorem{corollary}{Corollary}[section]
\newtheorem{conjecture}{Conjecture}[section]

\theoremstyle{definition}

\theoremstyle{remark}
\newtheorem{remark}{Remark}[section]

\numberwithin{equation}{section}
\makeatletter
\@namedef{subjclassname@2020}{%
  \textup{2020} Mathematics Subject Classification}
\makeatother
 \vspace{4mm}

\begin{document}
\hbox{To appear in Ramanujan J.}
\medskip

\title[On some determinants involving the tangent function]
{On some determinants involving the tangent function}
\author{Zhi-Wei Sun}
\address {Department of Mathematics, Nanjing
University, Nanjing 210093, People's Republic of China}
\email{zwsun@nju.edu.cn}

\keywords{Determinants, quadratic residues modulo primes, quadratic fields, the tangent function.
\newline \indent Supported by the Natural Science Foundation of China (grant no. 12371004).}

\begin{abstract}
Let $p$ be an odd prime and let $a,b\in\Bbb Z$ with $p\nmid ab$. In this paper we mainly evaluate
$$T_p^{(\delta)}(a,b,x):=\det\left[x+\tan\pi\frac{aj^2+bk^2}p\right]_{\delta\ls j,k\ls (p-1)/2}\ \ (\delta=0,1).$$
For example, in the case $p\eq3\pmod4$ we show that $T_p^{(1)}(a,b,0)=0$ and
$$T_p^{(0)}(a,b,x)=\begin{cases} 2^{(p-1)/2}p^{(p+1)/4}&\text{if}\ (\frac{ab}p)=1,
\\p^{(p+1)/4}&\text{if}\ (\frac{ab}p)=-1,\end{cases}$$
where $(\frac{\cdot}p)$ is the Legendre symbol. When $(\f{-ab}p)=-1$, we also evaluate the determinant $\det[x+\cot\pi\frac{aj^2+bk^2}p]_{1\ls j,k\ls(p-1)/2}.$
In addition, we pose several conjectures one of which states that
for any prime $p\eq3\pmod4$ there is an integer $x_p\eq1\pmod p$
such that $$\det\l[\sec2\pi\f{(j-k)^2}p\r]_{0\ls j,k\ls p-1}=-p^{(p+3)/2}x_p^2.$$
\medskip

\end{abstract}

\maketitle

\section{Introduction}

Let $p$ be an odd prime.
 It is well known that the numbers
$$0^2,\ 1^2,\ \ldots,\ \l(\f{p-1}2\r)^2$$
are pairwise incongruent modulo $p$. In \cite{S19a} the author investigated the determinants
$$S(d,p)=\det\l[\l(\f{j^2+dk^2}p\r)\r]_{1\ls j,k\ls(p-1)/2}$$
and
$$T(d,p)=\det\l[\l(\f{j^2+dk^2}p\r)\r]_{0\ls j,k\ls(p-1)/2},$$
where $d$ is an integer not divisible by $p$, and $(\f{\cdot}p)$ is the Legendre symbol.
In particular, Sun \cite{S19a} showed that if $(\f dp)=1$ then
$$\l(\f{-S(d,p)}p\r)=1\ \ \t{and}\ \ T(d,p)=\f{p-1}2S(d,p).$$

Inspired by the determinants $S(d,p)$ and $T(d,p)$ with $d\in\Z$ and $p\nmid d$,
and noting that the tangent function $\tan x$ has period $\pi$,
for $a,b\in\Z$ we introduce
\begin{equation}\label{1.1}T_p^{(0)}(a,b,x):=\det\l[x+\tan\pi\f{aj^2+bk^2}p\right]_{0\ls j,k\ls(p-1)/2}
\end{equation}
and
\begin{equation}\label{1.2}T_p^{(1)}(a,b,x):=\det\l[x+\tan\pi\f{aj^2+bk^2}p\right]_{1\ls j,k\ls(p-1)/2},
\end{equation}
and denote $T_p^{(0)}(a,b,0)$ and $T_p^{(1)}(a,b,0)$ by $T_p^{(0)}(a,b)$ and $T_p^{(1)}(a,b)$ respectively.
To study the novel determinants $T_p^{(0)}(a,b,x)$ and $T_p^{(1)}(a,b,x)$, we first find their values
by numerical experiments via {\tt Mathematica}, and then seek for detailed proofs via related known results involving roots of unity.

Now we present our main results.

\begin{theorem}\label{Th1.1} Let $p$ be an odd prime and let $a,b\in\Z$ with $p\nmid ab$.

{\rm (i)} Assume that $p\eq1\pmod4$. Then
\begin{equation}\label{1.3}T_p(a,b,-x)=-T_p^{(0)}(a,b,x),
\end{equation}
and in particular $T_p^{(0)}(a,b)=0.$
If $(\f {ab}p)=1$ and $b\eq ac^2\pmod p$ with $c\in\Z$, then
\begin{equation}\label{1.4}T_p^{(1)}(a,b,x)=\l(\f{2c}p\r)p^{(p-3)/4}\ve_p^{(\f ap)(2-(\f2p))h(p)},
\end{equation}
where $\ve_p$ and $h(p)$ are the fundamental unit and the class number of the real quadratic field $\Q(\sqrt p)$ respectively.
When $(\f{ab}p)=-1$, we have
\begin{equation}\label{1.5}T_p^{(1)}(a,b,x)=T_p^{(1)}(a,b)=\pm2^{(p-1)/2}p^{(p-3)/4}.
\end{equation}

{\rm (ii)} Suppose that $p\eq3\pmod 4$.
Then
\begin{equation}\label{1.6}T_p^{(1)}(a,b,-x)=-T_p^{(1)}(a,b,x),\end{equation}
and in particular $T_p^{(1)}(a,b)=0$.
Also,
\begin{equation}\label{1.7}T_p^{(0)}(a,b,x)=\begin{cases} 2^{(p-1)/2}p^{(p+1)/4}&\t{if}\ (\f{ab}p)=1,
\\p^{(p+1)/4}&\t{if}\ (\f{ab}p)=-1.\end{cases}
\end{equation}
\end{theorem}
\begin{remark}\label{Rem1.1} When $p$ is a prime with $p\eq1\pmod 4$, and $a$ and $b$ are integers with $(\f{ab}p)=-1$,
we are unable to determine the sign in \eqref{1.5}.
For any prime $p\eq3\pmod4$ and integers $a$ and $b$ with $p\nmid ab$, the identity \eqref{1.7}
looks surprising and interesting. We believe that Theorem \ref{1.1} has certain potential applications.
\end{remark}

\begin{theorem}\label{Th1.2} Let $n>1$ be an odd integer, and let $a$ and $b$ be integers with $\gcd(ab,n)=1$.
Then
\begin{equation}\label{1.8}\det\l[x+\tan\pi\f{aj+bk}n\r]_{0\ls j,k\ls n-1}+\det\l[-x+\tan\pi\f{aj+bk}n\r]_{0\ls j,k\ls n-1}=0
\end{equation}
and
\begin{equation}\label{1.9}\det\l[x+\tan\pi\f{aj+bk}n\r]_{1\ls j,k\ls n-1}=\l(\f{-ab}n\r)n^{n-2},
\end{equation}
where $(\frac{\cdot}n)$ is the Jacobi symbol.
\end{theorem}

For the cotangent function, we establish the following two theorems.

\begin{theorem}\label{Th1.3} Let $p>3$ be a prime, and let $a,b\in\Z$ with $(\f{-ab}p)=-1$. Then
\begin{equation}\label{1.10}\begin{aligned}&\ \det\l[x+\cot\pi\f{aj^2+bk^2}p\r]_{1\ls j,k\ls(p-1)/2}
\\=&\ \begin{cases} T_p^{(1)}(a,b)/(-p)^{(p-1)/4}=\pm 2^{(p-1)/2}/\sqrt p&\t{if}\ p\eq1\pmod4,
\\(-1)^{(h(-p)+1)/2}(\f ap)2^{(p-1)/2}/\sqrt p&\t{if}\ p\eq3\pmod4,
\end{cases}\end{aligned}\end{equation}
where $h(-p)$ is the class number of the imaginary quadratic field $\Q(\sqrt p\,\i)$ with $\i=\sqrt{-1}$.
\end{theorem}
\begin{remark}\label{Rem1.2} It is known that $2\nmid h(-p)$ for each prime $p\eq3\pmod4$.
In 1961 Mordell \cite{M61} even proved that
for any prime $p>3$ with $p\eq3\ (\mo\ 4)$ we have
$$\f{p-1}2!\eq(-1)^{(h(-p)+1)/2}\pmod p.$$
\end{remark}

\begin{theorem}\label{Th1.4} For any odd prime $p$, we have
$$D_p:=\det\l[\cot\pi\f{jk}p\r]_{1\ls j,k\ls(p-1)/2}\in\begin{cases}\Q&\t{if}\ p\eq1\pmod4,\\\sqrt{p}\ \Q&\t{if}\ p\eq3\pmod 4.
\end{cases}$$
\end{theorem}

We are going to provide several lemmas in the next section and then prove Theorem \ref{Th1.1} in Section 3.
Theorems \ref{Th1.2}--\ref{Th1.4} will be shown in Section 4.
In Section 5, we pose some conjectures on determinants involving the tangent function.

\section{Some Lemmas}

\begin{lemma}\label{Lem2.1} Let $A=[a_{jk}]_{0\ls j,k\ls n}$ be a matrix over a field.
Then
\begin{equation}\label{2.1}\det[x+a_{jk}]_{0\ls j,k\ls n}=\det A+x\det B,
\end{equation}
where $B=[b_{jk}]_{1\ls j,k\ls n}$ with $b_{jk}=a_{jk}-a_{j0}-a_{0k}+a_{00}$.
\end{lemma}
\Proof. As $(x+a_{jk})-(x+a_{0k})=a_{jk}-a_{0k}$ for all $0<j\ls n$ and $0\ls k\ls n$, we have
\begin{align*}\det[x+a_{jk}]_{0\ls j,k\ls n}=&\vmatrix x+a_{00}&x+a_{01}&x+a_{02}&\hdots&x+a_{0n}\\a_{10}-a_{00}&a_{11}-a_{01}&a_{12}-a_{02}&\hdots&a_{1n}-a_{0n}
\\\vdots&\vdots&\vdots&\ddots&\vdots\\a_{n0}-a_{00}&a_{n1}-a_{01}&a_{n2}-a_{02}&\hdots&a_{nn}-a_{0n}\endvmatrix
\\=&\vmatrix x&x&x&\hdots&x\\a_{10}-a_{00}&a_{11}-a_{01}&a_{12}-a_{02}&\hdots&a_{1n}-a_{0n}
\\\vdots&\vdots&\vdots&\ddots&\vdots\\a_{n0}-a_{00}&a_{n1}-a_{01}&a_{n2}-a_{02}&\hdots&a_{nn}-a_{0n}\endvmatrix
\\&+\vmatrix a_{00}&a_{01}&a_{02}&\hdots&a_{0n}\\a_{10}-a_{00}&a_{11}-a_{01}&a_{12}-a_{02}&\hdots&a_{1n}-a_{0n}
\\\vdots&\vdots&\vdots&\ddots&\vdots\\a_{n0}-a_{00}&a_{n1}-a_{01}&a_{n2}-a_{02}&\hdots&a_{nn}-a_{0n}
\endvmatrix,
\end{align*}
and hence $\det[x+a_{jk}]_{0\ls j,k\ls n}-\det A$ coincides with
$$\vmatrix x&0&\hdots&0\\a_{10}-a_{00}&a_{11}-a_{01}-(a_{10}-a_{00})&\hdots&a_{1n}-a_{0n}-(a_{10}-a_{00})
\\\vdots&\vdots&\ddots&\vdots\\a_{n0}-a_{00}&a_{n1}-a_{01}-(a_{n0}-a_{00})&\hdots
&a_{nn}-a_{0n}-(a_{n0}-a_{00})\endvmatrix=x\det B.$$
This concludes the proof of \eqref{2.1}. \qed

\begin{corollary}\label{Cor2.1} Let $m$ and $n$ be positive integers with $2\nmid n$.
Let $f:\Z\to \R$ be an odd function, where $\R$ is the field of real numbers. Then, for any integer $d$, the determinant
$$\det\left[x+f((j+d)^m-(k+d)^m)\right]_{0\ls j,k\ls n}$$
does not depend on $x$.
\end{corollary}
\Proof. Let
$$a_{jk}=f((j+d)^m-(k+d)^m)\quad\ \t{for}\ j,k=0,\ldots,n.$$
For $1\ls j,k\ls n$ set $b_{jk}=a_{jk}-a_{j0}-a_{0k}+a_{00}$. As $f$ is an odd function, we have
\begin{align*} b_{jk}&=f((j+d)^m-(k+d)^m)-f((j+d)^m-d^m)-f(d^m-(k+d)^m)
\\&=-f((k+d)^m-(j+d)^m)+f((k+d)^m-d^m)+f(d^m-(j+d)^m)=-b_{kj}.
\end{align*}
Thus
$$\det[b_{jk}]_{1\ls j,k\ls n}=(-1)^{n}\det[b_{kj}]_{1\ls j,k\ls n}=-\det[b_{jk}]_{1\ls j,k\ls n}$$
and hence $\det[b_{jk}]_{1\ls j,k\ls n}=0$. Applying Lemma 2.1, we immediately get the desired result. \qed

The following lemma is Frobenius' extension (cf. \cite{Frob} and \cite[(8)]{O}) of Cauchy's determinant identity (cf. \cite[(5.5)]{K05}).

\begin{lemma}\label{Lem2.2} We have
\begin{equation}\label{2.2}
\det\bg[z+\f1{x_j+y_k}\bg]_{0\ls j,k\ls n}=\f{\prod_{0\ls j<k\ls n}(x_k-x_j)(y_k-y_j)}{\prod_{j=0}^n\prod_{k=0}^n(x_j+y_k)}\(1+z\sum_{k=0}^n(x_k+y_k)\).
\end{equation}
\end{lemma}
\proof We present here an induction proof of \eqref{2.2} by using Cauchy's determinant identity
which is the special case $z=0$ of \eqref{2.2}.

In the case $n=0$, both sides of \eqref{2.2} coincide with $z+1/(x_0+y_0)$.

Now, let $n$ be a positive integer, and suppose that
\begin{equation}\label{2.2'}\det\bg[z+\f1{x_j+y_k}\bg]_{1\ls j,k\ls n}=\f{\prod_{1\ls j<k\ls n}(x_k-x_j)(y_k-y_j)}{\prod_{j=1}^n\prod_{k=1}^n(x_j+y_k)}\(1+z\sum_{k=1}^n(x_k+y_k)\).
\end{equation}
By Lemma 2.1 and \eqref{2.2} in the case $z=0$,
\begin{equation}\label{bjk}\det\bg[z+\f1{x_j+y_k}\bg]_{0\ls j,k\ls n}=\f{\prod_{0\ls j<k\ls n}(x_k-x_j)(y_k-y_j)}{\prod_{j=0}^n\prod_{k=0}^n(x_j+y_k)}+z\det[b_{jk}]_{1\ls j,k\ls n},
\end{equation}
where
$$b_{jk}=\f1{x_j+y_k}-\f1{x_j+y_0}-\f1{x_0+y_k}+\f1{x_0+y_0}=\f{(x_j-x_0)(y_k-y_0)(x_j+y_k+x_0+y_0)}
{(x_0+y_0)(x_j+y_0)(x_0+y_k)(x_j+y_k)}.$$
With the aid of \eqref{2.2'}, we have
\begin{align*}\det[b_{jk}]_{1\ls j,k\ls n}&=\prod_{j=1}^n\f{x_j-x_0}{x_j+y_0}
\times\prod_{k=1}^n\f{y_k-y_0}{y_k+x_0}\times\det\bg[\f1{x_0+y_0}+\f1{x_j+y_k}\bg]_{1\ls j,k\ls n}
\\&=\prod_{k=1}^n\f{(x_k-x_0)(y_k-y_0)}{(x_k+y_0)(y_k+x_0)}
\times\f{\prod_{1\ls j<k\ls n}(x_k-x_j)(y_k-y_j)}
{\prod_{j=1}^n\prod_{k=1}^n(x_j+y_k)}\l(1+\f{\sum_{k=1}^n(x_k+y_k)}{x_0+y_0}\r)
\\&=\f{\prod_{0\ls j<k\ls n}(x_k-x_j)(y_k-y_j)}{\prod_{j=0}^n\prod_{k=0}^n(x_j+y_k)}\sum_{k=0}^n(x_k+y_k).
\end{align*}
Combining this with \eqref{bjk}, we obtain the desired \eqref{2.2}. This concludes the proof. \qed

An analogue of Lemma \ref{Lem2.2} for Pfaffians can be found in Okada's paper \cite{O}.

\begin{lemma} [Huang and Pan \cite{P06}] \label{Lem2.3} Let $n>1$ be an odd integer, and let $c$ be any integer relatively prime to $n$. For each $j=1,\ldots,(n-1)/2$, let $\rho_c(j)$ be the unique $r\in\{1,\ldots,(n-1)/2\}$ with $cj$
congruent to $r$ or $-r$ modulo $n$. For the permutation $\rho_c$ on $\{1,\ldots,(n-1)/2\}$, its sign is given by
\begin{equation}\label{2.3}\sign(\rho_c)=\l(\f cn\r)^{(n+1)/2}.
\end{equation}
\end{lemma}

\begin{lemma}[Sun \cite{S19b}]\label{Lem2.4} Let $p$ be an odd prime. Let $\zeta=e^{2\pi \i/p}$ and $a\in\Z$ with $p\nmid a$.

{\rm (i)} If $p\eq1\pmod4$, then
\begin{equation}\label{2.4}\prod_{1\ls j<k\ls(p-1)/2}(\zeta^{aj^2}+\zeta^{ak^2})
=\pm\ve_p^{(\f ap)h(p)((\f2p)-1)/2}
\end{equation}
and
\begin{equation}\label{2.5}\prod_{1\ls j<k\ls(p-1)/2}(\zeta^{aj^2}-\zeta^{ak^2})^2
=(-1)^{(p-1)/4}p^{(p-3)/4}\ve_p^{(\f ap)h(p)}.
\end{equation}

{\rm (ii)} Suppose that $p\eq3\pmod4$. Then
\begin{equation}\label{2.6}\prod_{1\ls j<k\ls(p-1)/2}(\zeta^{aj^2}+\zeta^{ak^2})=1,
\end{equation}
and
\begin{equation}\label{2.7}\prod_{1\ls j<k\ls(p-1)/2}(\zeta^{aj^2}-\zeta^{ak^2})
= \begin{cases}(-p)^{(p-3)/8}&\t{if}\ p\eq3\pmod8,
\\(-1)^{(p+1)/8+(h(-p)-1)/2}(\f ap)p^{(p-3)/8}\i&\t{if}\ p\eq7\pmod8.
\end{cases}\end{equation}
Also,
\begin{equation}\label{2.8}\prod_{k=1}^{(p-1)/2}(1-\zeta^{ak^2})=(-1)^{(h(-p)+1)/2}\l(\f ap\r)\sqrt {p}\,\i.\end{equation}
\end{lemma}

\begin{lemma}\label{Lem2.5} Let $p$ be an odd prime, and let $a,b\in\Z$ with $(\f{-ab}p)=-1$. Then
\begin{equation}\label{2.9}\prod_{j=1}^{(p-1)/2}\prod_{k=1}^{(p-1)/2}\l(1-e^{2\pi\i(aj^2+bk^2)/p}\r)
= p^{(p-1)/4}\times\begin{cases}1&\t{if}\ p\eq1\pmod4,
\\(-1)^{(h(-p)-1)/2}(\f ap)\i&\t{if}\ p\eq3\pmod4.
\end{cases}\end{equation}
\end{lemma}
\Proof. For $m\in\Z$ set
\begin{align*} r(m):&=\l|\l\{(j,k):\ 1\ls j,k\ls\f{p-1}2\ \t{and}\ aj^2+bk^2\eq m\pmod p\r\}\r|
\\&=\l|\l\{1\ls x\ls p-1:\ \l(\f xp\r)=1\  \t{and}\ \l(\f{m-ax}p\r)=\l(\f bp\r)\r\}\r|.
\end{align*}
Note that $r(0)=0$ since $(\f{-ab}p)\not=1$.

Let $m\in\{1,\ldots,p-1\}$. Then
\begin{align*} r(m)=\ &\sum_{0<x<p\atop p\nmid ax-m}\f{(\f xp)+1}2\cdot\f{(\f{b(m-ax)}p)+1}2
\\=\ &\f14\sum_{x=1}^{p-1}\l(\l(\f{bx(m-ax)}p\r)+\l(\f xp\r)+\l(\f{b(m-ax)}p\r)+1\r)-\f{(\f{am}p)+1}4
\\=\ &\f14\sum_{x=0}^{p-1}\l(\f{-abx^2+bmx}p\r)+\f14\sum_{x=0}^{p-1}\l(\f xp\r)
+\f14\sum_{x=0}^{p-1}\l(\f{-abx+bm}p\r)
\\&-\f14\l(\f {bm}p\r)+\f{p-1}4-\f{(\f{am}p)+1}4.
\end{align*}
It is well known that for any $a_0,a_1,a_2\in\Z$ with $p\nmid a_0$ or $p\nmid a_1$ we have
\begin{equation}\label{2.10}\sum_{x=0}^{p-1}\l(\f{a_0x^2+a_1x+a_2}p\r)=\begin{cases}-(\f {a_0}p)&\t{if}\ p\nmid a_1^2-4a_0a_2,
\\(p-1)(\f {a_0}p)&\t{if}\ p\mid a_1^2-4a_0a_2.\end{cases}
\end{equation}
(See, e.g., \cite[p.\,58]{BEW}.) Therefore
$$r(m)=-\f14\l(\f{-ab}p\r)+\f{p-1}4-\f{(\f{am}p)+(\f{bm}p)+1}4=\f{p-1}4-\f{1-(\f{-1}p)}4\l(\f{am}p\r).$$

In view of the above,
\begin{align*}&\ \prod_{j=1}^{(p-1)/2}\prod_{k=1}^{(p-1)/2}\l(1-e^{2\pi\i(aj^2+bk^2)/p}\r)
\\=&\ \prod_{m=1}^{p-1}(1-e^{2\pi\i m/p})^{r(m)}
=\f{\prod_{m=1}^{p-1}(1-e^{2\pi\i m/p})^{(p-1+(\f ap)(1-(\f{-1}p)))/4}}
{\prod_{0<m<p\atop(\f mp)=1}(1-e^{2\pi\i m/p})^{(\f ap)(1-(\f{-1}p))/2}}.
\end{align*}
Clearly,
$$\prod_{m=1}^{p-1}(1-e^{2\pi\i m/p})=\lim_{x\to1}\f{x^p-1}{x-1}=p.$$
As \eqref{2.8} holds for $p\eq3\pmod4$, we have
$$\prod_{0<m<p\atop(\f mp)=1}(1-e^{2\pi\i m/p})^{(1-(\f{-1}p))/2}
=\begin{cases}1&\t{if}\ p\eq1\pmod4,\\(-1)^{(h(-p)+1)/2}\sqrt p\,\i&\t{if}\ p\eq3\pmod4.\end{cases}$$
Thus the desired \eqref{2.9} follows. \qed

\section{Proof of Theorem \ref{Th1.1}}

We can easily verify Theorem \ref{Th1.1} for $p=3$. Below we assume that $p>3$.
For convenience, we set $n=(p-1)/2$ and $\zeta=e^{2\pi\i/p}$.
Since
$$\sum_{k=0}^nk^2=\f{n(n+1)(2n+1)}6=\f{p^2-1}{24}p\eq0\pmod p,$$
we have
\begin{equation}\label{3.1}\prod_{k=0}^n\zeta^{k^2}=1.\end{equation}
As
$$\tan x=\f{2\sin x}{2\cos x}=\f{(e^{\i x}-e^{-\i x})/\i}{e^{\i x}+e^{-\i x}}
=\f{-\i(e^{2\i x}-1)}{e^{2\i x}+1}=-\i+\f{2\i}{e^{2\i x}+1},$$
we also have
\begin{equation}\label{3.2}\i+\tan\pi\f{ aj^2+bk^2}p=\f{2\i}{\zeta^{aj^2+bk^2}+1}\quad \t{for all}\ j,k=0,\ldots,n.\end{equation}

For each $\da\in\{0,1\}$ and integer $d\not\eq0\pmod p$, we claim that
\begin{equation}\label{3.3}T_p^{(\da)}(a,\pm ad^2,x)=\l(\f dp\r)^{n+1}T_p^{(\da)}(a,\pm a,x).
\end{equation}
We now explain this.
For $k=1,\ldots,n$ let $\rho_d(k)$ be the unique $r\in\{1,\ldots,n\}$ with $dk$ congruent to $r$ or $-r$ modulo $p$. In view of Lemma \ref{Lem2.3},
\begin{align*} T_p^{(1)}(a,\pm ad^2,x)&=\sum_{\tau\in S_n}\sign(\tau)\prod_{j=1}^n\l(x+\tan\pi\f{aj^2\pm a(d\tau(j))^2}p\r)
\\&=\sign(\rho_d)\sum_{\tau\in S_n}\sign(\rho_d\tau)\prod_{j=1}^n\l(x+\tan\pi\f{aj^2\pm a\rho_d(\tau(j))^2}p\r)
\\&=\l(\f dp\r)^{n+1}T_p^{(1)}(a,\pm a,x).
\end{align*}
If we extend the function $\rho_d$ by defining $\rho_d(0)=0$,
then the new $\rho_d$ is a permutation of $\{0,1,\ldots,n\}$ and its sign
is the same as the old one.
So, \eqref{3.3} also holds for $\da=0$.

\medskip
\noindent{\it Proof of the First Part of Theorem \ref{Th1.1}}. As $p\eq1\pmod4$, we have $n=(p-1)/2\eq0\pmod2$.
For $q=n!$ we have $q^2\eq-1\pmod p$ by Wilson's theorem, hence
\begin{align*} -T_p^{(0)}(a,b,x)&=\det\l[-x-\tan\pi\f{aj^2+bk^2}p\r]_{0\ls j,k\ls n}
\\&=\det\l[-x+\tan\pi\f{a(qj)^2+b(qk)^2}p\r]_{0\ls j,k\ls n}
=T_p^{(0)}(a,b,-x)\end{align*}
and thus $\det T_p^{(0)}(a,b)=0$.
\medskip

{\it Case} 1. $(\f{ab}p)=1$.

In this case, $b\eq ac^2\pmod p$ for some integer $c\not\eq0\pmod p$. Note that
$b\eq -a(qc)^2\pmod p$ and hence
$$T_p^{(1)}(a,b,x)=\l(\f {2c}p\r)T_p^{(1)}(a,-a,x)$$
by \eqref{3.3} and the equality $(\f qp)=(\f 2p)$ (cf. \cite[Lemma 2.3]{S19a}).

By Corollary \ref{Cor2.1},
\begin{align*}&\det\l[x+\tan\pi\f{aj^2-ak^2}p\r]_{1\ls j,k\ls n}
=\det\l[x+\tan\pi\f{a(j+1)^2-a(k+1)^2}p\r]_{0\ls j,k\ls n-1}
\end{align*}
does not depend on $x$. So, with the aid of \eqref{3.2}, we get
\begin{align*} T_p^{(1)}(a,-a,x)&=\det\l[\i+\tan\pi\f{aj^2-ak^2}p\r]_{1\ls j,k\ls n}
\\&=\det\l[\f{2\i}{e^{2\pi\i a(j^2-k^2)/p}+1}\r]_{1\ls j,k\ls n}
\\&=\prod_{k=1}^n(2\i\zeta^{ak^2})\times\det\l[\f1{\zeta^{aj^2}+\zeta^{ak^2}}\r]_{1\ls j,k\ls n}.
\end{align*}
(Recall that $\zeta=e^{2\pi\i/p}$.)
In light of Lemma \ref{Lem2.2},
$$\det\l[\f1{\zeta^{aj^2}+\zeta^{ak^2}}\r]_{1\ls j,k\ls n}=\f{\prod_{1\ls j<k\ls n}(\zeta^{ak^2}-\zeta^{aj^2})^2}
{\prod_{k=1}^n(\zeta^{ak^2}+\zeta^{ak^2})\times\prod_{1\ls j<k\ls n}(\zeta^{aj^2}+\zeta^{ak^2})^2}.$$
Therefore,
\begin{align*} T_p^{(1)}(a,-a,x)&=\i^{n}\prod_{1\ls j<k\ls n}\l(\f{\zeta^{ak^2}-\zeta^{aj^2}}{\zeta^{ak^2}+\zeta^{aj^2}}\r)^2
\\&=(-1)^{(p-1)/4}\f{\prod_{1\ls j<k\ls n}(\zeta^{ak^2}-\zeta^{aj^2})^2}{\prod_{1\ls j<k\ls n}(\zeta^{aj^2}+\zeta^{ak^2})^2}
=p^{(p-3)/4}\ve_p^{(\f ap)(2-(\f 2p))h(p)}
\end{align*}
with the aid of Lemma \ref{Lem2.4}(i).
\medskip

{\it Case}\ 2. $(\f{ab}p)=-1$.

Recall that $T_p^{(1)}(a,b)=\det[c_{jk}]_{1\ls j,k\ls n}$ with
$c_{jk}=\tan\pi(aj^2+bk^2)/p$. By Lemma \ref{Lem2.1},
\begin{equation}\label{3.4}T_p^{(1)}(a,b,x)=\det[x+c_{jk}]_{1\ls j,k\ls n}=T_p^{(1)}(a,b)+x\det[d_{jk}]_{1<j,k\ls n},\end{equation}
where $d_{jk}=c_{jk}-c_{j1}-c_{1k}+c_{11}$. In light of \eqref{3.2} and \eqref{3.4},
$$\det\l[\f{2\i}{\zeta^{aj^2+bk^2}+1}\r]_{1\ls j,k\ls n}=\det[\i+c_{jk}]_{1\ls j,k\ls n}
=T_p^{(1)}(a,b)+\i\det[d_{jk}]_{1<j,k\ls n},$$
and hence \eqref{1.5} is implied by
\begin{equation}\label{3.5}
D_p(a,b):=\det\l[\f{2\i}{\zeta^{aj^2+bk^2}+1}\r]_{1\ls j,k\ls n}=\pm2^{(p-1)/2}p^{(p-3)/4}.
\end{equation}
(Note that both $T_p^{(1)}(a,b)$ and $\det[d_{jk}]_{1<j,k\ls n}$ are real numbers.)

In view of Lemma \ref{Lem2.2} and \eqref{3.1},
\begin{align*} D_p(a,b)
&=\prod_{k=1}^n\l(\f{2\i}{\zeta^{bk^2}}\r)\times\det\l[\f1{\zeta^{aj^2}+\zeta^{-bk^2}}\r]_{1\ls j,k\ls n}
\\&=\f{(2\i)^{n}}{\prod_{k=1}^n\zeta^{bk^2}}\times\f{\prod_{1\ls j<k\ls n}(\zeta^{ak^2}-\zeta^{aj^2})(\zeta^{-bk^2}-\zeta^{-bj^2})}{\prod_{j=1}^n\prod_{k=1}^n(\zeta^{aj^2}+\zeta^{-bk^2})}
\\&=(-1)^{(p-1)/4}2^{(p-1)/2}\f{\prod_{1\ls j<k\ls n}(\zeta^{ak^2}-\zeta^{aj^2})(\zeta^{-bk^2}-\zeta^{-bj^2})}
{\prod_{j=1}^n\prod_{k=1}^n(\zeta^{aj^2+bk^2}+1)}.
\end{align*}
Note that
$$\prod_{j=1}^n\prod_{k=1}^n(\zeta^{aj^2+bk^2}+1)
=\prod_{j=1}^n\prod_{k=1}^n\f{1-\zeta^{2aj^2+2bk^2}}{1-\zeta^{aj^2+bk^2}}=1$$
by Lemma \ref{Lem2.5}. So
\begin{equation}\label{3.6}D_p(a,b)=(-1)^{(p-1)/4}2^{(p-1)/2}\prod_{1\ls j<k\ls n}(\zeta^{ak^2}-\zeta^{aj^2})(\zeta^{-bk^2}-\zeta^{-bj^2}).\end{equation}
Observe that
$$\prod_{1\ls j<k\ls n}(\zeta^{ak^2}-\zeta^{aj^2})^2(\zeta^{-bk^2}-\zeta^{-bj^2})^2=p^{(p-3)/2}\ve_p^{((\f ap)+(\f{-b}p))h(p)}=p^{(p-3)/2}$$
by \eqref{2.5}. Therefore \eqref{3.5} holds and hence so does \eqref{1.5}.

In view of the above, we have completed the proof of part (i) of Theorem \ref{Th1.1}. \qed

\medskip
\noindent{\it Proof of the Second Part of Theorem \ref{Th1.1}}.
As $p\eq3\pmod 4$, we have
$n=(p-1)/2\eq1\pmod 2.$

{\it Case} 1. $(\f{ab}p)=-1$.

In this case, $b\eq-ad^2\pmod p$ for some integer $d\not\eq0\pmod p$, and hence by \eqref{3.3} we have
$$T_p^{(0)}(a,b,x)=T_p^{(0)}(a,-a,x)\ \ \t{and}\ \ T_p^{(1)}(a,b,x)=T_p^{(1)}(a,-a,x).$$
As
$$T_p^{(1)}(a,-a,-x)=\det\l[-x+\tan\pi\f{ak^2-aj^2}p\r]_{1\ls j,k\ls n}=(-1)^nT_p^{(1)}(a,-a,x)=-T_p^{(1)}(a,-a,x),$$
we get $T_p^{(1)}(a,b,-x)=-T_p^{(1)}(a,b,x)$, and in particular $T_p^{(1)}(a,b)=0$.

To obtain the equality $T_p^{(0)}(a,b,x)=p^{(p+1)/4}$, we now determine $T_p^{(0)}(a,-a,x)$
which equals $T_p^{(0)}(a,b,x)$. In view of Corollary \ref{Cor2.1} and \eqref{3.2}, we have
\begin{align*} T_p^{(0)}(a,-a,x)&=\det\l[\i+\tan\pi\f{aj^2-ak^2}p\r]_{0\ls j,k\ls n}
\\&=\det\l[\f{2\i}{\zeta^{a(j^2-k^2)}+1}\r]_{0\ls j,k\ls n}
\\&=\prod_{k=0}^n(2\i\zeta^{ak^2})\times\det\l[\f1{\zeta^{aj^2}+\zeta^{ak^2}}\r]_{0\ls j,k\ls n}.
\end{align*}
By Lemma \ref{Lem2.2},
$$\det\l[\f1{\zeta^{aj^2}+\zeta^{ak^2}}\r]_{0\ls j,k\ls n}=\f{\prod_{0\ls j<k\ls n}(\zeta^{ak^2}-\zeta^{aj^2})^2}
{\prod_{k=0}^n(\zeta^{ak^2}+\zeta^{ak^2})\times\prod_{0\ls j<k\ls n}(\zeta^{aj^2}+\zeta^{ak^2})^2}.$$
Therefore,
\begin{align*} T_p^{(0)}(a,-a,x)&=\i^{n+1}\prod_{0\ls j<k\ls n}\l(\f{\zeta^{ak^2}-\zeta^{aj^2}}{\zeta^{ak^2}+\zeta^{aj^2}}\r)^2
\\&=(-1)^{(p+1)/4}\prod_{k=1}^n\l(\f{\zeta^{ak^2}-1}{\zeta^{ak^2}+1}\r)^2\times\f{\prod_{1\ls j<k\ls n}(\zeta^{ak^2}-\zeta^{aj^2})^2}{\prod_{1\ls j<k\ls n}(\zeta^{aj^2}+\zeta^{ak^2})^2}.
\end{align*}
By Lemma \ref{Lem2.4}(ii),
$$\prod_{k=1}^n(\zeta^{ak^2}-1)^2=-p\ \ \t{and}\
\ \prod_{k=1}^n(\zeta^{ak^2}+1)^2=\prod_{k=1}^n\f{(\zeta^{2ak^2}-1)^2}{(\zeta^{ak^2}-1)^2}=\f{-p}{-p}=1,$$
and $$\prod_{1\ls j<k\ls n}(\zeta^{ak^2}-\zeta^{aj^2})^2=(-p)^{(p-3)/4}\
\t{and}\ \prod_{1\ls j<k\ls n}(\zeta^{ak^2}+\zeta^{aj^2})^2=1.$$
Therefore
$$T_p^{(0)}(a,-a,x)=(-1)^{(p+1)/4}(-p)(-p)^{(p-3)/4}=p^{(p+1)/4}$$
as desired.
\medskip

{\it Case} 2. $(\f{ab}p)=1$.

In this case, $b\eq ac^2\pmod p$ for some $c\in\Z$ with $p\nmid c$, and hence by \eqref{3.3} we have
$T_p^{(0)}(a,b,x)=T_p^{(0)}(a,a,x)$ and $T_p^{(1)}(a,b,x)=T_p^{(1)}(a,a,x)$ since $n+1$ is even.

Clearly $T_p^{(0)}(a,a)=\det[a_{jk}]_{0\ls j,k\ls n}$ with $a_{jk}=\tan\pi(aj^2+ak^2)/p$. By Lemma \ref{Lem2.1},
\begin{equation}\label{3.7}T_p^{(0)}(a,a,x)=\det[x+a_{jk}]_{0\ls j,k\ls n}=T_p^{(0)}(a,a)+x\det[b_{jk}]_{1\ls j,k\ls n}
\end{equation}
where
$$b_{jk}:=a_{jk}-a_{j0}-a_{0k}+a_{00}=\tan\pi\f{aj^2+ak^2}p-\tan\pi\f{aj^2}p-\tan\pi\f{ak^2}p.$$
Using the well known identity
$$\tan(x_1+x_2)=\f{\tan x_1+\tan x_2}{1-\tan x_1\tan x_2},$$
we obtain
$$b_{jk}=\tan\pi\f{aj^2}p\times\tan\pi\f{ak^2}p\times\tan\pi\f{aj^2+ak^2}p$$
and hence
\begin{equation}\label{3.8}\det[b_{jk}]_{1\ls j,k\ls n}=T_p^{(1)}(a,a)\prod_{j=1}^n\tan^2\pi\f{aj^2}p.
\end{equation}
In view of \eqref{3.2}, \eqref{3.7} and \eqref{3.8},
\begin{align*}\det\l[\f{2\i}{\zeta^{a(j^2+k^2)}+1}\r]_{0\ls j,k\ls n}&=\det[\i+a_{jk}]_{0\ls j,k\ls n}
=T_p^{(0)}(a,a)+\i T_p^{(1)}(a,a)\prod_{j=1}^n\tan^2\pi\f{aj^2}p.
\end{align*}
Thus
\begin{equation}\label{3.9}T_p^{(0)}(a,a)=2^{(p-1)/2}p^{(p+1)/4},\ T_p^{(1)}(a,a)=0\ \t{and}\ \det[b_{jk}]_{1\ls j,k\ls n}=0
\end{equation}
if and only if
\begin{equation}\label{3.10}\det\l[\f{2\i}{\zeta^{a(j^2+k^2)}+1}\r]_{0\ls j,k\ls n}=2^{(p-1)/2}p^{(p+1)/4}.\end{equation}
With the aid of Lemma \ref{Lem2.2},
\begin{align*} \det\l[\f{2\i}{\zeta^{a(j^2+k^2)}+1}\r]_{0\ls j,k\ls n}
&=\prod_{k=0}^n\f{2\i}{\zeta^{ak^2}}\times\det\l[\f1{\zeta^{aj^2}+\zeta^{-ak^2}}\r]_{0\ls j,k\ls n}
\\&=\f{(2\i)^{n+1}}{\prod_{k=0}^n\zeta^{ak^2}}\times\f{\prod_{0\ls j<k\ls n}(\zeta^{ak^2}-\zeta^{aj^2})(\zeta^{-ak^2}-\zeta^{-aj^2})}{\prod_{j=0}^n\prod_{k=0}^n(\zeta^{aj^2}+\zeta^{-ak^2})}.
\end{align*}
This, together with \eqref{3.1}, yields
\begin{equation}\label{3.11}\det\l[\f{2\i}{\zeta^{a(j^2+k^2)}+1}\r]_{0\ls j,k\ls n}
 (-1)^{(p+1)/4}2^{(p+1)/2}
\f{\prod_{0\ls j<k\ls n}(\zeta^{ak^2}-\zeta^{aj^2})(\zeta^{-ak^2}-\zeta^{-aj^2})}{\prod_{j=0}^n\prod_{k=0}^n(\zeta^{a(j^2+k^2)}+1)}.
\end{equation}

By Lemma \ref{Lem2.4}(ii),
\begin{align*}&\ \prod_{0\ls j<k\ls n}(\zeta^{ak^2}-\zeta^{aj^2})(\zeta^{-ak^2}-\zeta^{-aj^2})
\\=&\ \prod_{k=1}^n(\zeta^{ak^2}-1)(\zeta^{-ak^2}-1)\times\prod_{1\ls j<k\ls n}(\zeta^{aj^2}-\zeta^{ak^2})(\zeta^{-aj^2}-\zeta^{-ak^2})
\\=&\ p\times p^{(p-3)/4}=p^{(p+1)/4}.
\end{align*}
In view of Lemma \ref{Lem2.4}(ii) and Lemma \ref{Lem2.5},
\begin{align*}\prod_{j=0}^n\prod_{k=0}^n(\zeta^{a(j^2+k^2)}+1)
&=(\zeta^0+1)\prod_{j=1}^n\l(\f{1-\zeta^{2aj^2}}{1-\zeta^{aj^2}}\r)^2
\times\prod_{j=1}^n\prod_{k=1}^n\f{1-\zeta^{2a(j^2+k^2)/p}}{1-\zeta^{a(j^2+k^2)/p}}
\\&=2\l(\f 2p\r)^2\l(\f 2p\r)=2(-1)^{(p+1)/4}.
\end{align*}
Combining these with \eqref{3.11}, we get \eqref{3.10} and hence \eqref{3.9} holds.
In view of \eqref{3.7} and \eqref{3.9}, we finally obtain that
$$T_p^{(0)}(a,b,x)=T_p^{(0)}(a,a,x)=2^{(p-1)/2}p^{(p+1)/4}.$$

By the above, we have finished the proof of part (ii) of Theorem \ref{Th1.1}. \qed

\section{Proofs of Theorems \ref{Th1.2}-\ref{Th1.4}}

The following lemma is Frobenius' extension (cf. \cite{BC15}) of the Zolotarev lemma \cite{Z}.

\begin{lemma}\label{Lem4.1} Let $n$ be a positive odd integer, and let $a\in\Z$ be relatively prime to $n$.
For $j=0,\ldots,n-1$, let $\lambda_a(j)$
be the least nonnegative residue of $aj$ modulo $n$. Then the permutation $\lambda_a$ of $\{0,\ldots,n-1\}$ has the sign $\sign(\lambda_a)=(\f an)$.
\end{lemma}

We also need another lemma.

\begin{lemma}\label{Lem4.2} Let $n>1$ be an odd number and let $a\in\Z$ with $\gcd(a,n)=1$. Then
\begin{equation}\label{4.1}\prod_{1\ls j<k\ls n-1}\l(e^{2\pi\i ak/n}-e^{2\pi\i aj/n}\r)^2=(-1)^{(n-1)/2}n^{n-2}.
\end{equation}
\end{lemma}
\Proof. Let $\zeta=e^{2\pi\i a/n}$. Clearly,
\begin{equation}\label{4.2}\prod_{r=1}^{n-1}(1-\zeta^r)=\lim_{x\to1}\f{x^n-1}{x-1}=n
\end{equation}
and hence
\begin{align*}(-1)^{\bi{n-1}2}\prod_{1\ls j<k\ls n-1}(\zeta^k-\zeta^j)^2
&=\prod_{j=1}^{n-1}\prod_{k=1\atop k\not=j}^{n-1}(\zeta^j-\zeta^{k})
=\prod_{j=1}^{n-1}\prod_{k=1\atop k\not=j}^{n-1}\zeta^j(1-\zeta^{k-j})
\\&=\prod_{j=1}^{n-1}\(\f{(\zeta^j)^{n-2}}{1-\zeta^{-j}}\prod_{k=0\atop k\not=j}^{n-1}(1-\zeta^{k-j})\)
\\&=\f{\zeta^{(n-1)\sum_{j=0}^{n-1}j}}{\prod_{j=1}^{n-1}(\zeta^j-1)}\prod_{j=1}^{n-1}\prod_{r=1}^{n-1}(1-\zeta^r)
=n^{n-2}.
\end{align*}
So \eqref{4.1} holds. \qed

\medskip
\noindent{\it Proof of Theorem \ref{Th1.2}}. In view of Lemma \ref{Lem4.1}, for each $\da=0,1$ we have
\begin{align*}\det\l[x+\tan\pi\f{aj+bk}n\r]_{\da\ls j,k\ls n-1}=\l(\f an\r)\det\l[x+\tan\pi\f{j+bk}n\r]_{\da\ls j,k\ls n-1}
=\l(\f {-ab}n\r)D_n^{(\da)}(x),\end{align*}
where
$$D_n^{(\da)}(x):=\det\l[x+\tan\pi\f{j-k}n\r]_{\da\ls j,k\ls n-1}.$$

Since
\begin{align*} D_n^{(0)}(-x)&=\det\l[-x+\tan\pi\f{k-j}n\r]_{0\ls j,k\ls n-1}=\det\l[-x-\tan\pi\f{j-k}n\r]_{0\ls j,k\ls n-1}
\\&=(-1)^n\det\l[x+\tan\pi\f{j-k}n\r]_{0\ls j,k\ls n-1}=-D_n^{(0)}(x),
\end{align*} we have
$$\det\l[-x+\tan\pi\f{aj+bk}n\r]_{0\ls j,k\ls n-1}=-\l(\f{-ab}n\r)D_n^{(0)}(x)
=-\det\l[x+\tan\pi\f{aj+bk}n\r]_{0\ls j,k\ls n-1}$$
and hence \eqref{1.8} holds.

Now it remains to show that $D_n^{(1)}(x)=n^{n-2}$. Write $\zeta=e^{2\pi\i/n}$.
Similar to \eqref{3.2}, we have
$$\i+\tan\pi\f{j-k}n=\f{2\i}{\zeta^{j-k}+1}\quad\ \t{for all}\ j,k=1,\ldots,n-1.$$
Thus
$$D_n^{(1)}(\i)=\det\l[\f{2\i}{\zeta^{j-k}+1}\r]_{1\ls j,k\ls n-1}
=\prod_{k=1}^{n-1}(2\i\zeta^k)\times\det\l[\f1{\zeta^j+\zeta^k}\r]_{1\ls j,k\ls n-1}.
$$
By Lemma \ref{Lem2.2},
\begin{align*}\det\l[\f1{\zeta^j+\zeta^k}\r]_{1\ls j,k\ls n-1}&=\f{\prod_{1\ls j<k\ls n-1}(\zeta^k-\zeta^j)^2}
{\prod_{j=1}^{n-1}\prod_{k=1}^{n-1}(\zeta^j+\zeta^k)}
\\&=\f{\prod_{1\ls j<k\ls n-1}(\zeta^k-\zeta^j)^2}
{\prod_{k=1}^{n-1}(2\zeta^k)\times\prod_{1\ls j<k\ls n-1}(\zeta^k+\zeta^j)^2}.
\end{align*}
Therefore
$$D_n^{(1)}(\i)=\i^{n-1}\prod_{1\ls j<k\ls n-1}\f{(\zeta^k-\zeta^j)^4}{(\zeta^{2k}-\zeta^{2j})^2}=(-1)^{(n-1)/2}\prod_{1\ls j<k\ls n-1}\f{(\zeta^k-\zeta^j)^4}{(\zeta^{2k}-\zeta^{2j})^2}.$$
Combining this with Lemma \ref{Lem4.2}, we immediately get
$D_n^{(1)}(\i)=(n^{n-2})^2/n^{n-2}=n^{n-2}$. By Lemma \ref{2.1},
$$D_n^{(1)}(x)=D_n^{(1)}(0)+rx$$
for certain real number $r$. As $D_n^{(1)}(\i)=n^{n-2}$, we have
$D_n^{(1)}(0)=n^{n-2}$ and $r=0$.
Thus $D_n^{(1)}(x)=D_n^{(1)}(0)=n^{n-2}$ as desired.

The proof of Theorem \ref{Th1.2} is now complete. \qed

\medskip\noindent{\it Proof of Theorem \ref{Th1.3}}.
For any nonzero real number $x\not\in\pi\Z$, we obviously have
$$\cot x=\f{\cos x}{\sin x}=\f{(e^{\i x}+e^{-\i x})/2}{(e^{\i x}-e^{-\i x})/(2\i)}=\i+\f{2\i}{e^{2\i x}-1}.$$
Thus
$$-\i+\cot\pi\f{aj^2+bk^2}p=\f{2\i}{\zeta^{aj^2+bk^2}-1}\quad \t{for all}\ j,k=1,\ldots,n,$$
where $n=(p-1)/2$ and $\zeta=e^{2\pi\i/p}$. Let
$$C(x)=\det\l[x+\cot\pi\f{aj^2+bk^2}p\r]_{1\ls j,k\ls n}.$$
Then
$$C(-\i)=\det\l[\f{2\i}{\zeta^{aj^2+bk^2}-1}\r]_{1\ls j,k\ls n}=\prod_{k=1}^n(2\i\zeta^{-bk^2})\times\det\l[\f1{\zeta^{aj^2}-\zeta^{-bk^2}}\r]_{1\ls j,k\ls n}.$$
By Lemma \ref{Lem2.2},
\begin{align*}\det\l[\f1{\zeta^{aj^2}-\zeta^{-bk^2}}\r]_{1\ls j,k\ls n}
&=\f{\prod_{1\ls j<k\ls n}(\zeta^{ak^2}-\zeta^{aj^2})(-\zeta^{-bk^2}-(-\zeta^{-bj^2}))}
{\prod_{j=1}^n\prod_{k=1}^n(\zeta^{aj^2}-\zeta^{-bk^2})}
\\&=(-1)^{\bi n2}\f{\prod_{1\ls j<k\ls n}(\zeta^{ak^2}-\zeta^{aj^2})(\zeta^{-bk^2}-\zeta^{-bj^2})}
{(\prod_{k=1}^n\zeta^{-bk^2})^n\prod_{j=1}^n\prod_{k=1}^n(\zeta^{aj^2+bk^2}-1)}.
\end{align*}
Note that $\prod_{k=1}^n\zeta^{k^2}=1$ by \eqref{3.1}. So
\begin{equation}\label{4.3}C(-\i)=(2\i)^n
\f{(-1)^{\bi n2}\prod_{1\ls j<k\ls n}(\zeta^{ak^2}-\zeta^{aj^2})(\zeta^{-bk^2}-\zeta^{-bj^2})}
{(-1)^n\prod_{j=1}^n\prod_{k=1}^n(1-\zeta^{aj^2+bk^2})}.
\end{equation}
\medskip

{\it Case} 1. $p\eq1\pmod4$, i.e., $2\mid n$.

In this case, $(\f{ab}p)=(\f{-ab}p)=-1$. By Lemma \ref{Lem2.5},
$$\prod_{j=1}^n\prod_{k=1}^n(1-\zeta^{aj^2+bk^2})=p^{(p-1)/4}.$$
Combining this with \eqref{4.3} and \eqref{3.6}, we get
$$C(-\i)=\f{2^{(p-1)/2}}{p^{(p-1)/4}}
\prod_{1\ls j<k\ls n}(\zeta^{ak^2}-\zeta^{aj^2})(\zeta^{-bk^2}-\zeta^{-bj^2})
=\f{D_p(a,b)}{(-p)^{(p-1)/4}},$$
where $D_p(a,b)$ is defined as in \eqref{3.5}.
Thus
$$C(-\i)=\f{D_p(a,b)}{(-p)^{(p-1)/4}}=\f{T_p^{(1)}(a,b)}{(-p)^{(p-1)/4}}
=\pm\f{2^{(p-1)/2}}{\sqrt p}$$
with the aid of \eqref{1.5}. By Lemma 2.1, $C(x)=C(0)+rx$ for certain real number $r$.
Since $C(-\i)$ is real, we have $r=0$ and hence
$$C(x)=C(-\i)=\f{T_p^{(1)}(a,b)}{(-p)^{(p-1)/4}}=\pm\f{2^{(p-1)/2}}{\sqrt p}.$$
\medskip

{\it Case} 2. $p\eq3\pmod4$, i.e., $2\nmid n$.

In light of \eqref{2.7},
$$\prod_{1\ls j<k\ls n}(\zeta^{ak^2}-\zeta^{aj^2})(\zeta^{-bk^2}-\zeta^{-bj^2})=p^{(p-3)/4}.$$
Combining this with \eqref{2.9} and \eqref{4.3},  we obtain
\begin{align*}C(-\i)=(2\i)^n(-1)^{\bi n2}\f{p^{(p-3)/4}}{(-1)^n(-1)^{(h(-p)-1)/2}(\f ap)p^{(p-1)/4}i}
=\f{2^n\i(\i^2)^{(n-1)/2}(-1)^{(n-1)/2}}{(-1)^{(h(-p)+1)/2}(\f ap)\sqrt p\,i}
\end{align*}
and hence
$$C(-\i)=(-1)^{(h(-p)+1)/2}\l(\f ap\r)\f{2^{(p-1)/2}}{\sqrt p}$$
is a real number. Combining this with Lemma 2.1, we get that
$$C(x)=C(-\i)=(-1)^{(h(-p)+1)/2}\l(\f ap\r)\f{2^{(p-1)/2}}{\sqrt p}.$$

 In view of the above, we have completed the proof of Theorem \ref{Th1.3}. \qed
\medskip

 The following lemma is a well known result on quadratic Gauss sums (cf. \cite[pp.\,70-76]{IR}).

 \begin{lemma}\label{Lem4.3} Let $p$ be an odd prime. Then, for any integer $a\not\eq0\pmod p$, we have
 $$\sum_{x=0}^{p-1}e^{2\pi\i ax^2/p}=\l(\f ap\r)\sum_{t=0}^{p-1}\l(\f tp\r)e^{2\pi\i t/p}
 =\l(\f ap\r)\sqrt{(-1)^{(p-1)/2}p}.$$
 \end{lemma}

 Let $p$ be an odd prime, and let $\zeta=e^{2\pi\i/p}$. For $a,b\in\Z$ with $p\nmid ab$,
 Lemmas \ref{Lem2.2} and \ref{Lem4.3} are helpful to evaluate $\det[z+1/(\zeta^{aj^2}+\zeta^{bk^2})]_{1\ls j,k\ls(p-1)/2}$. However, we actually only need the case $z=0$ in our previous proofs of Theorems \ref{Th1.1}--\ref{Th1.3}.

 \medskip\noindent{\it Proof of Theorem \ref{Th1.4}}.  The Galois group $\text{Gal}(\Q(e^{2\pi\i/p})/\Q)$ consists of those automorphisms $\sigma_a$ $(1\le a\le p-1)$ with $\sigma_a(e^{2\pi\i/p})=e^{2\pi\i a/p}$. For any integer $x\not\eq0\pmod p$, we have
 $$\cot\pi\f xp=\i\f{e^{\pi\i x/p}+e^{-\pi\i x/p}}{e^{\pi\i x/p}-e^{-\pi\i x/p}}=\i\f{e^{2\pi\i x/p}+1}{e^{2\pi\i x/p}-1}.$$
 It follows that
 $$\f{D_p}{\i^{(p-1)/2}}=\det\l[\f{e^{2\pi\i jk/p}+1}{e^{2\pi\i jk/p}-1}\r]_{1\ls j,k\ls(p-1)/2}.$$

Let $a\in\{1,\ldots,p-1\}$. By the last equality,
\begin{align*}\sigma_a\l(\f{D_p}{\i^{(p-1)/2}}\r)&=\sigma_a\(\sum_{\tau\in S_{(p-1)/2}}\sign(\tau)\prod_{j=1}^{(p-1)/2}\f{e^{2\pi\i j\tau(j)/p}+1}{e^{2\pi\i j\tau(j)/p}-1}\)
\\&=\sum_{\tau\in S_{(p-1)/2}}\sign(\tau)\prod_{j=1}^{(p-1)/2}\f{e^{2\pi\i aj\tau(j)/p}+1}{e^{2\pi\i aj\tau(j)/p}-1}
\\&=\det\l[\f{e^{2\pi\i ajk/p}+1}{e^{2\pi\i ajk/p}-1}\r]_{1\ls j,k\ls(p-1)/2}
=\frac1{\i^{(p-1)/2}}\det\left[\cot\pi\frac{ajk}p\right]_{1\le j,k\le(p-1)/2}.
\end{align*}
By Gauss' Lemma (see, e.g., \cite[p.\,52]{IR}),
$$\left(\frac ap\right)=(-1)^{|\{1\le j\le(p-1)/2:\ \{aj/p\}>1/2\}|},$$
where $\{x\}$ denotes the fractional part of a real number $x$.
Therefore,
\begin{align*}\sigma_a\l(\f{D_p}{\i^{(p-1)/2}}\r)=&\frac{(\frac ap)}{\i^{(p-1)/2}}\det\left[\cot\pi\frac{\rho_a(j)k}p\right]_{1\le j,k\le(p-1)/2},
\end{align*}
where $\rho_a(j)$ is the unique $r\in\{1,\ldots,(p-1)/2\}$ with $aj\equiv \pm r\pmod p$.
Combining this with Lemma \ref{Lem2.3}, we deduce that
\begin{align*}\sigma_a\l(\f{D_p}{\i^{(p-1)/2}}\r)&=\frac{(\frac ap)}{\i^{(p-1)/2}}\l(\f ap\r)^{(p+1)/2}\det\left[\cot\pi\frac{jk}p\right]_{1\le j,k\le(p-1)/2}
\\&=\l(\f ap\r)^{(p-1)/2}\f{D_p}{\i^{(p-1)/2}}.
\end{align*}
If $p\eq1\pmod4$, then $\i^{(p-1)/2}=(-1)^{(p-1)/4}\in\Q$ and hence
$\sigma_a(D_p)=D_p$.
When $p\eq3\pmod4$, by Lemma \ref{Lem4.3} we have $\sqrt{-p}\in \Q(e^{2\pi\i/p})$
and
$$\sigma_a(\sqrt{-p})=\sigma_a\(\sum_{x=0}^{p-1}e^{2\pi\i x^2/p}\)
=\sum_{x=0}^{p-1}e^{2\pi\i ax^2/p}=\left(\frac ap\right)\sqrt{-p},$$
therefore
\begin{align*}\sigma_a\left((-1)^{(p+1)/4}\frac {D_p}{\sqrt{p}}\right)&=\sigma_a\left(\frac {D_p}{\i(-1)^{(p-3)/4}\sqrt{-p}}\right)=\f{\sigma_a(D_p/\i^{(p-1)/2})}{\sigma_a(\sqrt{-p})}
\\&=\frac {(\f ap)D_p/\i^{(p-1)/2}}{(\f ap)\sqrt{-p}}=(-1)^{(p+1)/4}\f{D_p}{\sqrt p}
\end{align*}
and hence $\sigma_a(D_p/\sqrt p)=D_p/\sqrt p$.

By the above, if $p\eq1\pmod4$, then $\sigma(D_p)=D_p$ for all $\sigma\in \t{Gal}(\Q(e^{2\pi\i/p}))$, and hence $D_p\in\Q$ by Galois theory. Similarly,
when $p\eq3\pmod4$ we have $D_p/\sqrt p\in\Q$.

The proof of Theorem \ref{Th1.4} is now complete. \qed

\section{Some open conjectures}

\begin{conjecture}\label{Conj5.1} Let $p$ be any odd prime. Then
\begin{equation}\label{5.1}\l(\f{-2}p\r)\f{\det\l[\cot\pi{jk}/p\r]_{1\ls j,k\ls (p-1)/2}}{2^{(p-3)/2}p^{(p-5)/4}}\in\{1,2,3,\ldots\},
\end{equation}
and this number is divisible by $h(-p)$ if $p\eq3\pmod4$.
\end{conjecture}
\begin{remark}\label{Rem5.1} By Theorem \ref{Th1.4}, for any odd prime $p$ we have
$$\det\left[\cot\pi\f{jk}p\right]_{1\ls j,k\ls(p-1)/2}\in p^{(p-1)/4}\Q.$$
\end{remark}

\begin{conjecture}\label{Conj5.2} Let $n$ be a positive integer.

{\rm (i)} The number
\begin{equation}\label{5.2}s_n:=(2n+1)^{-n/2}\det\l[\tan\pi\f{jk}{2n+1}\r]_{1\ls j,k\ls n}
\end{equation}
is always an integer.

{\rm (ii)} We have
\begin{equation}\label{5.3}\det\l[\tan^2\pi\f{jk}{2n+1}\r]_{1\ls j,k\ls n}\in(2n+1)^{(n+1)/2}4^{n-1}\Z.
\end{equation}
\end{conjecture}
\begin{remark}\label{Rem5.2} Via {\tt Mathematica} we find that
\begin{gather*} s_1=1,\ s_2=-2,\ s_3=s_4=4,\ s_5=48,\ s_6=-160,
\\ s_7=32,\ s_8=2176,\ s_9=6912,\ s_{10}=0,\ s_{11}=273408.
\end{gather*}
Let $t_n$ denote the $n$th term of the sequence \cite[A277445]{I16}, which is the determinant
of a matrix $T(n)=[t_{jk}]_{1\ls j,k\ls n}$ with entries among $0,\pm1$
such that
$$2\sum_{k=1}^nt_{jk}\sin\f{\pi k}{2n+1}=\tan\f{\pi j}{2n+1}\quad \ \t{for all}\ j=1,\ldots,n.$$
We guess that
$s_n=-t_n$ if $n\eq3\pmod4$, and $s_n=t_n$ otherwise.
\end{remark}

\begin{conjecture}\label{Conj5.3} For any odd integer $n>1$, we have
\begin{equation}\label{5.4}\det\l[\tan^2\pi\f{j+k}n\r]_{1\ls j,k\ls n-1}\in n^{n-2}\Z.
\end{equation}
\end{conjecture}
\begin{remark}\label{Rem5.3} We are able to prove that $\det[\tan^2\pi\f{j-k}n]_{1\ls j,k\ls n-1}\in\Z$ for any odd integer $n>1$.
\end{remark}

\begin{conjecture}\label{Conj5.4} Let $p\eq3\pmod4$ be a prime, and let $a,b\in\Z$ with $p\nmid ab$. Then
\begin{equation}\label{5.5}\det\l[\tan^2\pi\f{aj^2+bk^2}p\r]_{1\ls j,k\ls(p-1)/2}\in p^{(p-3)/4}\Z
\end{equation}
and
\begin{equation}\label{5.6}\det\l[\tan^2\pi\f{aj^2+bk^2}p\r]_{0\ls j,k\ls(p-1)/2}\in p^{(p+1)/4}\Z.\end{equation}
If $(\f{ab}p)=1$, then
\begin{equation}\label{5.7}\det\l[\cot^2\pi\f{aj^2+bk^2}p\r]_{1\ls j,k\ls(p-1)/2}\in \f{2^{p-3}}p\Z.
\end{equation}
\end{conjecture}

Let $p\eq1\pmod4$ be a prime, and let $a,b\in\Z$ with $p\nmid ab$. Choose $q\in\Z$ with $q^2\eq-1\pmod p$. Then
\begin{align*}&\ \det\l[\l(\f{a(qj)^2+b(qk)^2}p\r)\tan\pi\f{a(qj)^2+b(qk)^2}p\r]_{0\ls j,k\ls(p-1)/2}
\\=&\ (-1)^{(p+1)/2}\det\l[\l(\f{aj^2+bk^2}p\r)\tan\pi\f{aj^2+bk^2}p\r]_{0\ls j,k\ls(p-1)/2}
\end{align*}
and hence
\begin{equation}\label{5.8}\det\l[\l(\f{aj^2+bk^2}p\r)\tan\pi\f{aj^2+bk^2}p\r]_{0\ls j,k\ls(p-1)/2}=0.\end{equation}

\begin{conjecture}\label{Conj5.5} Let $p\eq3\pmod4$ be a prime and let $a,b\in\Z$ with $p\nmid ab$. Then
\begin{equation}\label{5.9}\det\l[\l(\f{aj^2+bk^2}p\r)\tan\pi\f{aj^2+bk^2}p\r]_{0\ls j,k\ls(p-1)/2}\in p\Z.\end{equation}
If $(\f{ab}p)=1$, then
\begin{equation}\label{5.10}\sqrt p\det\l[\l(\f{aj^2+bk^2}p\r)\cot\pi\f{aj^2+bk^2}p\r]_{1\ls j,k\ls(p-1)/2}\in \Z.\end{equation}
\end{conjecture}
\begin{remark}\label{Rem5.4} For any prime $p\eq3\pmod4$, set
$$a_p^{\pm}:=\f1p\det\l[\l(\f{j^2\pm k^2}p\r)\tan\pi\f{j^2\pm k^2}p\r]_{0\ls j,k\ls(p-1)/2}.$$
Via {\tt Mathematica} we find that
\begin{gather*} a_3^+=a_3^-=-1,\ a_7^+=60,\ a_7^-=3,\ a_{11}^+=2^6\times3^3,\ a_{11}^-=-373,
\\a_{19}^+=2^{12}\times3\times5^2\times 7\times11\times17
\ \t{and}\ a_{19}^-=-5\times7\times89\times3803.
\end{gather*}
\end{remark}

\begin{conjecture}\label{Conj5.6} Let $p$ be an odd prime.

{\rm (i)} Define
$$S(p):=\det\left[\sec 2\pi\f{jk}p\right]_{0\ls j,k\ls(p-1)/2}.$$
If $p\eq1\pmod4$, then $S(p)=0$. When $p\eq3\pmod4$, the number
$$\f{S(p)}{2^{(p-3)/2}(-p)^{(p+1)/4}}$$ is a positive odd integer.

{\rm (ii)} We have
$$c_p:=\f1{2^{(p-1)/2}p^{(p-5)/4}}\det\left[\csc2\pi\f{jk}p\right]_{1\ls j,k\ls(p-1)/2}\in\Z.$$
Moreover, $c_p=1$ if $p\eq3\pmod8$, and $c_p=0$ if $p\eq7\pmod 8$.
\end{conjecture}
\begin{remark}\label{Rem5.5} By the way we prove Theorem \ref{Th1.4}, we can show that
$S(p)/p^{(p+1)/4}\in\Q$ and $c_p\in\Q$
for any odd prime $p$.
In 2019 the author \cite{S19c} conjectured that
\begin{equation}\label{5.11}\f1{2n}\det\l[\cos\pi\f{jk}n\r]_{0\ls j,k\ls n}=\det\l[\cos\pi\f{jk}n\r]_{1\ls j,k\ls n}=(-1)^{\lfloor\f{n+1}2\rfloor}\f {n^{(n-1)/2}}{2^{(n-1)/2}}
\end{equation}
for every positive integer $n$, this was later confirmed by Petrov (cf. the answer in \cite{S19c}).
\end{remark}

\begin{conjecture}\label{Conj5.7} For any prime $p\eq3\pmod4$, there is an integer $x_p\eq1\pmod p$
such that
\begin{equation}\label{5.12}\det\l[\sec2\pi\f{(j-k)^2}p\r]_{0\ls j,k\ls p-1}=-p^{(p+3)/2}x_p^2.
\end{equation}
\end{conjecture}
\begin{remark} For $p=3,7,11$, we may take $x_p=1$ in \eqref{5.12}. For each prime $p\eq3\pmod4$,
the author \cite{MOsin} conjectured in 2021 that
$$\det\l[\sin2\pi\f{(j-k)^2}p\r]_{1\ls j,k\ls p-1}=-\f{p^{(p-1)/2}}{2^{p-1}},$$
which was later confirmed by Kalmynin (cf. the answer in \cite{MOsin}).
\end{remark}


\medskip


\begin{thebibliography}{99}


\bibitem{BC15} A. Brunyate and P. L. Clark, {\it Extending the Zolotarev-Frobenius approach to quadratic reciprocity}, Ramanujan J. {\bf 37} (2015), 25--50.

\bibitem{BEW} B. C. Berndt, R. J. Evans and K. S. Williams,
Gauss and Jacobi Sums, John Wiley \& Sons, 1998.

\bibitem{Frob} G. Frobenius, {\it \"Uber die elliptischen Funktionen zweiter Art}, J. Reine Angew. Math. {\bf 93} (1882), 53--68.

\bibitem{P06} C. Huang and H. Pan, {\it A remark on Zolotarev's theorem}, Colloq. Math. {\bf 171} (2023), 159--166.

\bibitem{I16} D. V. Ingerman, Sequence A277445 in OEIS,  {\tt http://oeis.org/A277445}, 2016.

\bibitem{IR} K. Ireland and M. Rosen, A Classical
Introduction to Modern Number Theory, Graduate Texts in
Math. 84, 2nd ed., Springer, New York, 1990.


\bibitem{K05} C. Krattenthaler, {\it Advanced determinant calculus: a complement},
Linear Algebra Appl. {\bf 411} (2005), 68--166.

\bibitem{M61} L. J. Mordell, {\it The congruence $((p-1)/2)!\eq\pm1\ (\mo\ p)$},
Amer. Math. Monthly {\bf 68} (1961), 145--146.

\bibitem{O} S. Okada, {\it An elliptic generalization of Schur's Pfaffian identity}, Adv. Math.
{\bf 204} (2006), 530--538.

\bibitem{S19a} Z.-W. Sun, {\it On some determinants with Legendre symbol entries},
Finite Fields Appl. {\bf 56} (2019), 285--307.

\bibitem{S19b} Z.-W. Sun, {\it Quadratic residues and related permutations and identities},
Finite Fields Appl. {\bf 59} (2019), 246--283.

\bibitem{S19c} Z.-W. Sun, {\it A surprising identity: $\det[\cos\pi jk/n]_{1\ls j,k\ls n}=(-1)^{\lfloor (n+1)/2\rfloor}(n/2)^{(n-1)/2}$},
    Question 321098 at MathOverflow (with an answer by Fedor Petrov), January 17, 2019.
Available from {\tt https://mathoverflow.net/questions/321098}

\bibitem{MOsin} Z.-W. Sun, {\it Is it true that $\det[\sin2\pi(j-k)^2/p]_{1\ls j,k\ls p-1}
=-p^{(p-1)/2}/2^{p-1}$?} Question 395143 at MathOverflow (with an answer by Alexander Kalmynin), June 12, 2021.
Available from {\tt https://mathoverflow.net/questions/395143}

\bibitem{Z} G. Zolotarev, {\it Nouvelle d\'emonstration de la loi de r\'eciprocit\'e de Legendre},
Nouvelles Ann. Math. {\bf 11} (1872), 354--362.

\end{thebibliography}
\end{document}